\documentclass[12pt]{article}
\usepackage[dvips]{graphicx}
\usepackage{t1enc,amssymb,amsbsy,amsthm,amsmath,amsfonts}

\renewcommand{\Re}{{\rm I\kern-0.16em R}}

\def\@begintheorem#1#2{\trivlist \item[\hskip \labelsep{\bf #1\ #2}]}
\def\@opargbegintheorem#1#2#3{\trivlist
      \item[\hskip \labelsep{\bf #1\ #2\ (#3)}]}

\newtheorem{proposition}{Proposition}[section] 

\newtheorem{thm}[proposition]{Theorem}
\newtheorem{lemma}[proposition]{Lemma}
\newtheorem{corollary}[proposition]{Corollary}
\newtheorem{example}[proposition]{Example}
\newtheorem{remark}[proposition]{Remark}

\def\P{{\bf P}}
\def\R{{\bf R}}
\def\R{{\bf R}}

\def\E{{\bf E}}

\def\cB{{\cal B}}

\def\cM{{\cal M}}

\def\al{\alpha}

\def\la{\lambda}

\numberwithin{equation}{section}

\begin{document}

\author{Ernesto Mordecki
\\{\small Universidad de la Rep\'ublica,}
\\{\small Cento de Matem\'atica,}
\\{\small Igu\'a 4225,}
\\{\small 11400 Montevideo, Uruguay}
\\{\small email: mordecki@cmat.edu.uy}
\and
Paavo Salminen
\\{\small \AA bo Akademi,}
\\{\small Mathematical Department,}
\\{\small F\"anriksgatan 3 B,}
\\{\small FIN-20500 \AA bo, Finland,} 
\\{\small email: phsalmin@abo.fi}
}
\title{Optimal stopping of Hunt \\ and L\'evy processes}
\maketitle
\begin{abstract}
The optimal stopping problem for a Hunt processes on $\R$ is considered via 
the representation theory of excessive functions. In particular, we focus on infinite
horizon (or perpetual) problems with one-sided structure, that is, there exists
a point $x^*$ such that the stopping region is of the form
$[x^*,+\infty)$. Corresponding results for two-sided problems are also
indicated. 
The main result is a spectral representation of the value function
in terms of the Green kernel of the process.
Specializing in L\'evy processes, 
we obtain, by applying the Wiener-Hopf factorization, a general
representation of the value function
in terms of the maximum of the L\'evy process.
To illustrate the results, an explicit expression for  
the Green kernel of Brownian motion with
exponential jumps is computed and some optimal stopping problems for Poisson
process with positive exponential jumps and negative drift are solved.
\end{abstract} 


\newpage
\section{Introduction}
\label{sec0}
Consider an optimal stopping problem for a real-valued Markov process
$X=\{X_t\colon t\geq 0\}$ 
with reward function $g$ and discount rate $r\geq 0$. 
Denote by $V$ the value function of the problem, and by $\tau^*$ the optimal
stopping time.
In this paper we analyze this optimal stopping problem 
departing from three main sources: (i) the characterization of
the value function $V$ as the least excessive majorant of the
reward function $g$, due to Snell  \cite{snell53} for discrete martingales
and to Dynkin \cite{dynkin63} for continuous time Markov processes; 
(ii) the representation of excessive functions  as integrals of the
Green kernel of the process, as exposed
in Kunita and Watanabe \cite{kunitawatanabe65} and
Dynkin \cite{dynkin69}, and exploited by Salminen \cite{salminen85}
in the framework of optimal stopping for diffusions; and (iii)
recent results expressing the solution of
some optimal stopping problems for L\'evy processes and random walks
in terms of the maximum of the process, 
see Darling et. al.\cite{darlingliggettaylor72}, 
Mordecki \cite{mordecki02},
Boyarchenko and Levendorskij \cite{boyarchenkolevendorskij02},
Novikov and Shiryaev \cite{novikovshiryaev04} and 
Kyprianou and Surya \cite{kyprianousurya05}.
For papers on optimal stopping of L\'evy processes  
using other methods, see, e.g., McKean  \cite{mckean67}, Gerber and Shiu
 \cite{gerbershiu98}, Chan \cite{chan00} and Kou and Wang  \cite{kouwang04}. 

We then try to understand the structure of the solution of
the optimal stopping problem in a regular enough framework
of Markov processes, precisely the class of Hunt processes, 
concluding that finding the solution of such a problem is equivalent
to finding the representation 
of the value function in terms of the Green kernel.
The Radon  measure that appears in
this representation is called the \emph{spectral measure} corresponding to
the excessive function $V$, and furthermore, it results that
the support of this spectral measure is
the stopping region for the problem. 
This is our main result, presented in section 3.

Let us specialize to L\'evy processes. 
Firstly, observe that in the case $r>0$
the Green kernel is proportional to the distribution
of the process stopped at an exponential time with parameter $r$,
independent of the process. 
Secondly, relying on the Wiener-Hopf factorization for the L\'evy process, 
we express this random variable in the distributational sense as the sum of two independent
random variables, the first one having the distribution of the supremum
of the process up to the exponential time and the second one 
with the distribution of the infimum of the process in
the same random interval. Simple calculations taking into account
this fact, and the one-sided structure of the
solution of the optimal stopping problem, gives 
a representation of the solution of an optimal
stopping problem in terms of the maximum of
the L\'evy process -- a result that has been
obtained earlier in several particular cases. This analysis
is carried out in section 4. 

The rest of the paper is as follows.
In section 2 the framework of Hunt processes
in which we are working is described.
Section 5 consists of two subsections. In the first one we 
illustrate the made assumptions concerning the
Hunt processes and L\'evy processes by studying Browian motion with
exponential jumps. In the second one an optimal stopping problem for a
compound Poisson process with negative drift and positive exponential
jumps and the reward functon $g(x)=(x^+)^\gamma,\ \gamma\geq 1,$ is analyzed.
Our interest in this particular reward function was arised by
Alexander Novikov's talk in the Symposium on Optimal Stopping with
Applications held in Manchester 22.-- 27.1.2006 \cite{novikov06}
where the optimal stopping problem for the same reward functions
and general L\'evy processes were considered.

\section{Preliminaries on Hunt processes}
\label{sec1}

Let $X=\{X_t\}$ be a transient Hunt process taking values
in $\R$, where we omit $t\geq 0$ in the notation,
as all the considered processes are indexed in the same set.
In particular, $X$ is a strong Markov process, quasi left
continuous on $[0,+\infty)$ and  the sample paths of $X$ are right
continuous with left limits (see Kunita and
Watanabe \cite{kunitawatanabe65} and Blumenthal and Getoor
\cite{blumenthalgetoor68} p. 45). The notations $\P_x$ and $\E_x$ are
used for the probabilty measure and the expectation operator, respectively, associated
with $X$ when $X_0=x.$ The resolvent
$\{G_r\,:\,r\geq 0\}$ of $X$ is defined via
\begin{equation}\label{resolvent}
 G_r(x,A):=\int_0^\infty {\rm e}^{-r t}\, \P_x(X_t\in A)\, dt, 
\end{equation}
where $x\in\R$ and $A$ is a Borel subset of $\R$. We assume also that
there exists a dual resolvent  $\{\widehat G_r\,:\,r\geq
0\}$ with respect to some $\sigma$-finite (duality) measure $m,$ that is, 
for all  $f,g\in{\cal B}_o$ and $r\geq 0$ it holds
$$
\int m(dx)\, f(x)\, G_r g(x)=\int m(dx)\, \widehat G_r f(x)\,g(x),
$$
where ${\cal B}_o$ denotes the set of measurable bounded functions with compact support. 
Moreover, it is assumed that  $\{\widehat G_r\,:\,r\geq 0\}$ is a
resolvent of a transient Hunt process $\widehat X$ taking values in $\R.$
Finally, we impose Hypothesis (B) from \cite{kunitawatanabe65}
p. 498:
\begin{description}
\item{(${\rm h}_1$)}\quad $G_r(x,dy)=G_r(x,y)\,m(dy),$
\item{(${\rm h}_2$)}\quad $\sup_{x\in A} G_0(x,B)<\infty$ for all compact $A$ and $B,$
\item{(${\rm h}_3$)}\quad $x\mapsto \widehat G_r f(x), \,f\in {\cal B}_o,$ is continuous and finite.
\end{description}
The assumption that the dual process is a
Hunt process implies that $\widehat G$ is regular (see 
\cite{kunitawatanabe65} p. 494).

We remark that when $X$ is a L\'evy process a dual resolvent always exists
as the resolvent of the dual process
$\widehat{X}=\{-X_t\}$. Hereby
the Lebesgue measure serves as the duality measure (see section II.1 in \cite{bertoin96}).

Under these assumptions it can be proved that the function $G_r$
given in (${\rm h}_1$) constitutes a potential kernel (often called a Green kernel) of
exponent $r$ associated with the pair $(X,\widehat X).$ This means
that for each given $r\geq 0$ the function $(x,y)\mapsto G_r(x,y)$ is
jointly measurable and has the properties    
\begin{description}
\item{(${\rm p}_1$)}\quad  $\widehat G_r(y,dx)=G_r(x,y)\,m(dx),$
\item{(${\rm p}_2$)}\quad $x\mapsto G_r(x,y)$ is $r$-excessive for $X$ for each $y,$
\item{(${\rm p}_3$)}\quad $y\mapsto G_r(x,y)$ is $r$-excessive for $\widehat X$ for each $x.$
\end{description}
Recall that a non-negative measurable function $f:\R\mapsto[0,+\infty]$ is called
$r$-excessive for $X$ if for all $x\in\R$
\begin{description}
\item{(${\rm e}_1$)}\quad ${\rm e}^{-rt}\,\E_x(f(X_t))\leq f(x)$ for all $t\geq 0,$ 
\item{(${\rm e}_2$)}\quad ${\rm e}^{-rt}\,\E_x(f(X_t))\to f(x)$ as $t\to 0.$
\end{description}
Notice that $r$-excessive functions of $X$ are 0-excessive for the
process $X^o$ obtained from $X$ by exponential killing with rate $r$.  

From the assumption (${\rm h}_1$) that the resolvent kernel of $X$ is absolutely
continuous it follows 
that $r$-excessive functions are lower semi-continuous. 

The Riesz decomposition of excessive functions is of key importance
in our approach to optimal stopping. We state the decomposition
relying on \cite{kunitawatanabe65} Theorem 2 p. 505 and Proposition
13.1 p. 523. Indeed, it holds, under the made assumptions,  
that each $r$-excessive function $u$ locally
integrable with respect to the duality measure $m$ can be decomposed
uniquely in the form
\begin{equation}
\label{riesz} 
u(x)=\int_\R G_r(x,y)\, \sigma_u(dy) + h_r(x),
\end{equation}
where $h_r$ is an $r$-harmonic function and $\sigma_u$ is a Radon
measure on $\R.$ We remark that  
the assumption that the dual process is a
Hunt process implies that also the spectral measure $\sigma_u$ is unique 
(see \cite{kunitawatanabe65} Proposition 7.11 p. 503).

Conversily (see  \cite{kunitawatanabe65} Proposition
7.6 p. 501), given a Radon measure $\sigma$ on
$(\R,\cB )$ the function $v$ defined via
$$
v(x):= \int_\R G_r(x,y)\, \sigma(dy)
$$
is an $r$-excessive function (in fact, a potential).

An $r$-excessive function $u$ is said to $r$-harmonic on a Borel
subset $A$ of $\R$ if for all open subsets $A_c$ of $A$ with compact closure 
\begin{equation}
\label{harm} 
u(x)=\E_x\left( {\rm e}^{-rH(A_c)}\,u(X_{H(A_c)})\right),
\end{equation}
 where
$$ 
H(A_c):=\inf\{t\,:\, X_t\not\in A_c\}.
$$
In our case (see \cite{kunitawatanabe65} Proposition 6.2 p. 499) it holds
that for each fixed $y\in\R$ the function
$x\mapsto G_r(x,y)$ is $r$-harmonic on $\R\setminus\{y\}.$
From (\ref{riesz}) it follows that if there
exists an open set $A$ such that the representing measure does not
charge $A,$  then the function $u$ is
$r$-harmonic on $A$ i.e., 
\begin{equation}
\label{harm2} 
\sigma_u(A)=0 \quad \Rightarrow \quad u\ {\rm is\ }r{\rm
  -harmonic\ on\ } A. 
\end{equation}
In fact, (\ref{harm2}) is an equivalence under mild assumptions; 
for this see Dynkin \cite{dynkin69} Theorem 12.1.

\section{Optimal stopping of Hunt processes}
\label{sec3}

Let $X$ be a Hunt process satisfying the assumptions made in Section 2. 
For simplicity we consider non-negative continuous reward functions. 
Then  the reward function $g$ has the smallest excessive majorant
$V$  and 
\begin{equation}
\label{dyn}
V(x)=\sup_{\tau\in{\cM}}\E_x({\rm e}^{-r\tau}g(X_\tau)),
\end{equation}
where ${\cM}$ denotes the set of all stopping times $\tau$ with
respect to the natural filtration generated by $X.$ 
In case $\tau=+\infty$ we define
$$
{\rm e}^{-r\tau}g(X_{\tau}):=
\limsup_{t\to\infty}{\rm e}^{-r t}g(X_t).
$$ 
The result \eqref{dyn} can be found, for instance in Shiryayev \cite{shiryayev78} 
(Lemma 3 p. 118 and Theorem 1 p. 124) and holds for more general standard Markov
processes, and for almost-Borel and ${\cal C}_0$-lower semicontinuous reward functions.
We can express this result by saying that in an optimal
stopping problem the value function and the smallest excessive
majorant of the reward function coincide.

From  (\ref{dyn}) and the Riesz decomposition (\ref{riesz}) we
conclude that 
the problem of finding the value function
is equivalent to the problem of finding the representing measure
of the smallest excessive majorant (up to harmonic functions). 
Furthermore, based on (\ref{harm2}),
it is seen, roughly speaking,  that the continuation region, that is,
the region where it
is not optimal to stop, is the ``biggest'' set
not charged by the representing measure $\sigma_V$ of $V.$
In short, the representing measure gives the value function
via (\ref{riesz}) and the stopping region is by (\ref{harm2}) the support of
the representing measure.  
In the following result we use the preceeding considerations
in order to express the solution of a particular
type of optimal stopping problems

\begin{thm}
\label{prop1}
Consider a Hunt process  $\{X_t\}$ satisfying the assumptions made in
Section 2, a non-negative continuous reward function $g$,
and a discount rate $r\geq 0$ such that 
\begin{equation}
\label{r1}
\E_x(\sup_{t\geq 0}e^{-rt}g(X_t))<\infty.
\end{equation}
Assume that there exists a Radon 
measure $\sigma$ 
with support on the set $[x^*,\infty)$ 
such that the
function
\begin{equation}
\label{excessive}
V(x):=\int_{[x^*,\infty)}G_r(x,y)\sigma(dy)
\end{equation}
satisfies the following conditions:
\begin{itemize}
\item[\rm (a)]
$V$ is  continuous,
\item[\rm (b)]
$V(x)\to 0$ when $x\to-\infty$.
\item[\rm (c)]
$V(x)=g(x)\ \text{when $x\geq x^*$}$,
\item[\rm (d)]
$V(x)\geq g(x)\ \text{when $x<x^*$}$.
\end{itemize}
Let 
\begin{equation}\label{tau-star}
\tau^*=\inf\{t\geq 0\colon X_t\geq x^*\}.
\end{equation}
Then $\tau^*$ is an optimal stopping time and $V$ is the value
function of the optimal stopping problem for $\{X_t\}$ with the reward function $g,$
in other words,
\begin{equation*}
V(x)= \sup_{\tau\in{\cal M}}\E_x\left({\rm e}^{-r\tau} g(X_{\tau})\right)=\E_x\left( {\rm e}^{-r\tau^*}g(X_{\tau^*})\right), 
\quad x\in\R.
\end{equation*}
\end{thm}
\begin{proof}
Since $V$  is an $r$-excessive function (see the discussion after the
Riesz decomposition (\ref{riesz})) and, from conditions (c) and (d),
a majorant of $g$, it follows
by Dynkin's characterization of the value function
as the least excessive majorant, that
\begin{equation}\label{equal}
V(x)\geq \sup_{\tau\in{\cal M}}\E_x\left({\rm e}^{-r\tau} g(X_{\tau})\right)
\end{equation}
In order to conclude the proof, we establish the equality in \eqref{equal}.
Indeed, consider for each $n\geq 1$ the stopping time
\begin{equation*}
\tau_n=\inf\{t\geq 0\colon X_t\not\in (-n, x^*-1/n)\}.
\end{equation*}
For $\omega\in\{\tau^*<\infty\}$ define $\bar{\tau}=\lim_{n\to\infty}\tau_n$.
We have
$$
\tau_1\leq\tau_2\leq\cdots\leq\bar{\tau}\leq\tau^*.
$$
For $n$ large enough, $X_{\tau_n}\geq x^*-1/n$,
and, as the process is quasi-left continous, 
$\lim_ {n\to\infty}X_{\tau_n}=X_{\bar{\tau}}$, and, hence, 
$X_{\bar{\tau}}\geq x^*$. This give us that $\bar{\tau}=\tau^*$ a.s.

As $V$ is $r$-excessive, the
sequence $\{{\rm e}^{-r\tau_n}V(X_{\tau_n})\}$ is a nonnegative supermartingale,
and, in consequence, 
it converges a.s. to a random variable. Because $X_{\tau_n}\to X_{\tau^*}$ a.s.,
and $V$ is continuous, we identify the limit as ${\rm e}^{-r\tau^*}V(X_{\tau^*})$.
From assumptions (a) and (b) it follows that 
$$
C_V:=\sup_ {x\leq x^*} V(x)<\infty,
$$
Furthermore, as 
\begin{multline*}
e^{-r\tau_n}V(X_{\tau_n})
        =e^{-r\tau_n}V(X_{\tau_n}){\bf 1}_{\{{\tau_n}<{\tau_*}\}}
+e^{-r\tau^*}g(X_{\tau^*}){\bf 1}_{\{{\tau_n}={\tau_*}\}}\\
\leq C_V+\sup_{t\geq 0}e^{-rt}g(X_t)
\end{multline*}
we obtain, in view of condition \eqref{r1}, 
using the Lebesgue dominated convergence theorem
$$
\E_x\left({\rm e}^{-r\tau_n}\,V(X_{\tau_n})\right)\downarrow \E_x\left({\rm e}^{-r\tau^*}\, V(X_{\tau^*})\right)
\quad {\rm as\ } n\to\infty.
$$
Furthermore, as the representing measure $\sigma$ does not charge 
the open set $(-\infty,x^*)$, the function $V$ is harmonic on that
set (cf. (\ref{harm2})),
and, as $\tau_n$ are exit times from the open sets $(-n,x^*-1/n)$,
we conclude that 
\begin{equation*}
V(x)=\E\left({\rm e}^{-r\tau_n}\, V(X_{\tau_n})\right)
\downarrow
\E\left({\rm e}^{-r\tau^*}\, g(X_{\tau^*})\right),
\end{equation*}
and the proof is complete.
\end{proof}
Under the additional assumption \eqref{visit}, 
valid in many particular cases, we characterize now the optimal
threshold $x^*$ as a solution, with a useful uniqueness property, of an equation derived from 
(\ref{excessive}).  Remember that by the definition of the support
of a Radon measure, we have $\sigma(A)>0$ for all open subsets $A$ of $[x^*,\infty)$
(see for instance page 215 in Folland \cite{folland99}). 
\begin{corollary} 
\label{cor1}
Let $\{X_t\}, g, V, \sigma,$ and $x^*$ be as in Theorem
\ref{prop1}. Assume that
\begin{equation}  
\label{visit}
G_r(x,B):=\int_0^\infty {\rm e}^{-rt}\P_x(X_t\in B)dt\, >\,0
\end{equation}  
for all  $x$ and open subsets $B$ of $\R.$ Then the equation
\begin{equation}  
\label{visit1}
g(x)=\int_{[x,\infty)}G_r(x,y)\sigma(dy)
\end{equation}  
has no solution bigger than $x^*.$ 
\end{corollary}
\begin{proof} Clearly, since $g(x^*)=V(x^*)$ it is immediate from
  (\ref{excessive}) that $x^*$ is a solution of (\ref{visit1}). 
Let now $x_o>x^*$ be another solution of (\ref{visit1}), i.e.,
$$
g(x_o)=\int_{[x_o,\infty)}G_r(x_o,y)\sigma(dy).
$$ 
From (\ref{excessive}) we have
$$
g(x_o)=V(x_o)=\int_{[x^*,\infty)}G_r(x_o,y)\sigma(dy). 
$$
Consequently,
\begin{equation}  
\label{visit2}
\int_{[x^*,x_o)}G_r(x_o,y)\sigma(dy)=0. 
\end{equation}  
But the function $y\mapsto G_r(x_o,y)$ is lower semi-continuous, and,
hence, the set $\{y: G_r(x_o,y)>0\}$ is open. From (\ref{visit2}) it
follows that $G_r(x_o,\cdot)\equiv 0$ on $(x^*,x_o),$ but this
violates (\ref{visit}) and the claim is proved.
\end{proof}
The presented method works similarly 
when the stopping region is not a half line, i.e. when
the problem is not a ``one-sided'' problem.
The form of the optimal stopping time \eqref{tau-star} 
appears very often in several optimal stopping problems,
in particular in mathematical finance,
where this sort of the problems are sometimes named \emph{call-like}
perpetual problems or options (see e.g. \cite{boyarchenkolevendorskij02}).
Furthermore, as exposed in section
\ref{section:levy}, one sided problems in the context of
L\'evy process admits a representation in terms of
the maximum of the L\'evy process.

Minor modifications in the proof of Theorem \ref{prop1}
give the following result,
that can be considered as a ``two-sided'' optimal stopping problem.
\begin{thm}
\label{prop2}
Consider a Hunt process  $\{X_t\}$ satisfying the assumptions made in
Section 2 and a non-negative continuous reward function $g$,
and a discount rate $r\geq 0$, such that 
condition \eqref{r1} hold.
Assume that there exists a positive Radon 
measure $\sigma$ with support on the set $S=(-\infty,x_*]\cup [x^*,\infty)$ 
such that the
function
\begin{equation}
\label{excessive2}
V(x):=\int_SG_r(x,y)\sigma(dy)=
\int_{(-\infty,x_*]}G_r(x,y)\sigma(dy)+
\int_{[x^*,\infty)}G_r(x,y)\sigma(dy)
\end{equation}
satisfies the following conditions:
\begin{itemize}
\item[\rm (a)]
$V$ is  continuous,
\item[\rm (b)]
$V(x)=g(x)\ \text{when $x\in S$}$,
\item[\rm (c)]
$V(x)\geq g(x)\ \text{when $x\notin S$}$.
\end{itemize}
Let 
\begin{equation*}
\tau^*=\inf\{t\geq 0\colon X_t\in S\}.
\end{equation*}
Then $\tau^*$ is an optimal stopping time and $V$ is the value
function of the optimal stopping problem for $\{X_t\}$ with the reward function $g$,
in other words,
\begin{equation*}
V(x)= 
\sup_{\tau\in{\cal M}}\E_x\left({\rm e}^{-r\tau} g(X_{\tau})\right)
=\E_x\left( {\rm e}^{-r\tau^*}g(X_{\tau^*})\right), 
\quad x\in\R.
\end{equation*}
\end{thm}

\section{Optimal stopping and maxima for L\'evy processes}
\label{section:levy}
As we have mentioned, in several papers explicit solutions
to optimal stopping problems for general random walks
or L\'evy process, and some particular reward functions
can be expressed in terms of the maximum of the process,
killed at a constant rate $r$, the discount rate of the problem.
The pioneer results in this direction are contained in the paper 
of Darling, Ligget and Taylor \cite{darlingliggettaylor72}, 
were solutions to optimal stopping problems
for rewards $g(x)=({\rm e}^x-1)^+$ and $g(x)=x^+$
are obtained in the whole class of random walks, in terms
of the maximum of the random walk. These results
are generalized for L\'evy process by Mordecki 
in
\cite{mordecki02}
and 
\cite{mordecki02a},
where it is also observed that similar
results hold for solutions
of optimal stopping problems
in terms of the infimum of the process
for the payoff $g(x)=(1-{\rm e}^x)^+$. 
Based on the technique of factorization of
pseudo-differential operators, Boyarchenko and
Levendorskij (see
\cite{boyarchenkolevendorskij02}
and the references therein)
obtain similar results, 
in a subclass of L\'evy processes, called  
regular and of exponential type (RLPE),
with the important feature that their results
are not based on particular properties
of the reward function $g$, and, hence, hold true
in a certain class of rewards functions.
Results in \cite{boyarchenkolevendorskij02} suggest that any optimal
stopping problem 
for a L\'evy process
with an increasing
payoff can be expressed 
in terms of the maximum of the process.
More recently, Novikov and Shiryaev
\cite{novikovshiryaev04}
obtained the solution 
of the optimal stopping problem for a general
random walk,
in terms of the
maximum, when the reward is
$g(x)=(x^+)^n$, (and also when $g(x)=1-{\rm e}^{-x^+}$).
The respective generalization of 
this problem to the framework of L\'evy processes
has been carried out by Kyprianou and Surya \cite{kyprianousurya05}.

\subsection{L\'evy processes}

Let $X=\{X_t\}$ be a Hunt process with
independent and stationary increments, i.e. a L\'evy process.
We denote $\E=\E_0$ and $\P=\P_0$.


If $v\in\R$, L\'evy-Khinchine formula states $\E({\rm e}^{iv X_t})= {\rm e}^{t\psi(iv)}$,
where, for complex $z=iv$ the characteristic exponent of the process is
\begin{equation}
\psi(z)=az+\frac 12b^2z^2+
\int_{\R}\big({\rm e}^{zx}-1-zh(x)\big)\Pi(dx).
\label{eq:char-exp}
\end{equation}
Here the truncation function $h(x)=x{\bf 1}_{\{|x|\leq 1\}}$ is fixed, 
and the parameters characterizing the law of the process are:
the drift $a$, an arbitrary real number;
the standard deviation of the Gaussian part of the process $b\geq 0$;
and the L\'evy jump measure $\Pi$, a non negative measure, defined on ${\R}\setminus\{0\}$
such that 
$\int (1\wedge x^2)\Pi(dx)<+\infty$.


Denote by $\tau(r)$ an exponential random variable with parameter
$r>0$,   independent   of the process   $X$.  
A key role in this section is played by 
the following random variables:
\begin{equation}
M_r=\sup_{0\leq  t<\tau(r)}X_t\qquad\mbox{and}
\qquad  
I_r=\inf_{0\leq t<\tau(r)}X_t 
\label{max-min}
\end{equation}
called the \emph{supremum}  and the \emph{infimum} of  the  process,  
respectively,
killed at  rate $r$.

A relevant instrument to study these random variables is the
Wiener-Hopf-Rogozin factorization, obtained by Rogozin \cite{rogozin66}, 
that states
\begin{equation}\label{eq:wh}
\frac{r}{r -\psi(iv)}=
\E({\rm e}^{ivM_r})
\E({\rm e}^{ivI_r})
\end{equation}

In our first result we give some simple sufficient conditions in order to
hypothesis \eqref{r1} to hold.
\begin{lemma}
\label{lem1}
Assume that a non-negative function $g$ satisfies
\begin{equation}
\label{exp}
g(x)\leq A_0+A_1{\rm e}^{\alpha x},
\end{equation}
for nonnegative constants $A_0, A_1, \alpha$.
Assume furthermore 
that 
\begin{equation}
\label{exponentialmoment}
\E\left({\rm e}^{\alpha X_1}\right)<{\rm e}^r.
\end{equation}
Then, condition \eqref{r1} holds.
\end{lemma}
\begin{proof} Let us first verify that, for $r\geq 0$, 
the following three statements are equivalent:
\begin{itemize}
\item[\rm (a)] $\E\left({\rm e}^{\alpha X_1}\right) < {\rm e}^r$.
\item[\rm (b)] $\E\left({\rm e}^{\alpha M_r}\right) < \infty$.
\item[\rm (c)] $\E\left(\sup_{t\geq 0}\left({\rm e}^{\alpha X_t-rt}\right)\right) < \infty$.
\end{itemize}
First, (a)$\Leftrightarrow$(b) is Lemma 1 in \cite{mordecki02}.
The equivalence (a)$\Leftrightarrow$(c) is a particular case of  (a)$\Leftrightarrow$(b),
when considering the L\'evy process $\{\alpha X_t-rt\}$,
the first constant equal to 1, and the second, the discount rate equal to $0$.
Now
\begin{multline*}
\E_x(\sup_{t\geq 0}{\rm e}^{-rt}g(X_t))
\leq
\E\left(
\sup_{t\geq 0}{\rm e}^{-rt}\left(A_0+A_1{\rm e}^{\alpha(x+X_t)}
\right)
\right)
\\
\leq 
A_0+A_1{\rm e}^{\alpha x}\E(\sup_{t\geq 0}{\rm e}^{(\alpha X_t-rt)})<\infty
\end{multline*}
as condition (a)$\Rightarrow$(c).
\end{proof}
\begin{remark}
Condition \eqref{exp} is relatively natural in our context.
For instance, if the function is increasing, 
and \emph{submultiplicative}
(as defined in section 25 in \cite{sato99})
it automatically satisfies our exponential growth condition \eqref{exp}.
Nevertheless, the submultiplicative property does not seem
to be appropiate for optimal stopping problems, as $g(x)=x^+$ is
not submultiplicative.  
Furthermore, condition \eqref{exponentialmoment} is optimal in the following
sense: For the reward function $g(x)=({\rm e}^x-1)^+$, if
$\E\left({\rm e}^{X_1}\right)={\rm e}^r$, then 
condition \eqref{r1} does not hold, based on {\rm (a)}$\Leftrightarrow${\rm (c)}.
\end{remark}
Our next result represents the
value function of the optimal stopping
problem for a L\'evy process in terms of the 
maximum of the process and is a consequence of Theorem \ref{prop1}.
\begin{proposition}
\label{corollary:onesided}
Assume that the conditions of Theorem  \ref{prop1} hold, and, 
furthermore, that $\{X_t\}$ is in fact a L\'evy process.
Then, there exists a function $Q\colon [x^*,\infty)\to\R$ such
that the value function $V$ in  \eqref{excessive} 
satisfies
\begin{equation*}
V(x)=\E_x\left(Q(M_r)\,;\,M_r\geq x^*\right),
\qquad
x\leq x^*.
\end{equation*}
\end{proposition}
\begin{proof}
The key ingredient of the proof is formula \eqref{eq:wh},
that can be also written as
\begin{equation}
\label{whr}
X_{\tau(r)}=M_r+\tilde{I}_r
\end{equation}
where $M_r$ and $\tilde{I}_r=X_{\tau(r)}-M_r$ are 
\emph{independent}
random variables, $M_r$ given in \eqref{max-min},
and $\tilde{I}_r$ with the same distribution as $I_r$ in \eqref{max-min}.

From the definition of the Green kernel \eqref{resolvent},
it is clear that
$$
rG_r(x,dy)=\P_x(X_{\tau(r)}\in dy),
$$
and, in view of \eqref{whr}, assuming that $M_r$ and $I_r$ have
respective densities $f_M$ and $f_I$ (only for simplicity of exposition), 
we obtain that
\begin{equation}
rG_r(x,y)=\begin{cases}
\int_{-\infty}^{y-x}f_I(t)f_M(y-x-t)dt,&\text{if $y-x<0$},\\
{}&\text{}\\
\int_{y-x}^{\infty}f_M(t)f_I(y-x-t)dt,&\text{if $y-x>0$}.\\
\end{cases}
\end{equation}
If we plugg in this formula for the Green kernel
in \eqref{excessive}, when $x<x^*$, and, in consequence,
with $y>x$, we obtain
\begin{align*}
V(x)&=\int_{x^*}^{\infty}G_r(x,y)\sigma(dy)\\
&=
r^{-1}\int_{x^*}^{\infty}\left[
\int_{y-x}^{\infty}f_M(t)f_I(y-x-t)dt\right]
\sigma(dy)\\
&=r^{-1}\int_{x^*-x}^{\infty}f_M(t)
\left[
\int_{x^*}^{x+t}f_I(y-x-t)\sigma(dy)
\right]dt\\
&=
\int_{x^*-x}^{\infty}f_M(t)Q(x+t)dt=\E_x\left(Q(M_r)\,;\,M_r\geq x^*\right),
\end{align*}
where, for $z\geq x^*$, we denote
\begin{equation}\label{q}
Q(z)=r^{-1}\int_{x^*}^{z}f_I(y-z)\sigma(dy).
\end{equation}
This concludes the proof.
\end{proof}
The following results uses Theorem \ref{prop2}
to provide a representation of the value function
in terms of both the supremum and the infimum of
the L\'evy process.
\begin{proposition}
Assume that the conditions of Theorem  \ref{prop2} hold, and 
that $\{X_t\}$ is a L\'evy process.
Then, there exist two functions 
$$
Q_*\colon (-\infty,x_*]\to\R,
\qquad 
Q^*\colon [x^*,\infty)\to\R
$$
such that the value function $V$ in \eqref{excessive} 
satisfies 
\begin{equation*}
V(x)=\E_x\left(Q_*(I_r)\,;\,I_r\leq x_*\right)+
\E_x\left(Q^*(M_r)\,;\,M_r\geq x^*\right),
\qquad
x_*\leq x\leq x^*.
\end{equation*}
\end{proposition}
\begin{proof} 
The proof consist in rewriting each summand
in \eqref{excessive} in terms of the maximum and infimum of
the process, respectively.
The second identity has been obtained in Proposition \ref{corollary:onesided},
and states (with $Q^*$ instead of $Q$), that
$$
\int_{[x^*,\infty)}G_r(x,y)\sigma(dy)=
\E_x\left(Q^*(M_r)\,;\,M_r\geq x^*\right),
$$
where $Q^*$ is defined in \eqref{q}.
The first one is obtained from this last equality 
considering the dual L\'evy process $\widehat{X}$, as follows:
\begin{multline*}
\int_{(-\infty,x_*]}G_r(x,y)\sigma(dy)=\int_{[-x_*,\infty)}\widehat{G}_r(-x,y)\sigma(-dy)\\
=\widehat{\E}_{-x}\left(\widehat{Q}^*(\widehat{M}_r)\,;\,\widehat{M}_r\geq -x_*\right)
=\E_{x}\left(Q_*(I_r)\,;\,I_r\leq x_*\right),
\end{multline*}
where
$$
Q_*(z)=\widehat{Q}^*(-z)=
r^{-1}\int_{-x_*}^{-z}f_{\widehat{I}}(y+z)\sigma(-dy)=
r^{-1}\int_z^{x_*}f_M(y-z)\sigma(dy),
$$
and the proof is complete.
\end{proof}

\section{A case study}
\label{sec4}
\subsection{Brownian motion with exponential jumps}
\label{sec41}
Here we illustrate  the assumptions made in Section \ref{sec1} and, in
particular, the concept of Green kernel by 
taking $X$ to be a Brownian motion with drift and compounded
with two-sided exponentially distributed jumps. 

To introduce $X$, consider
a standard Wiener process $W=\{W_t\colon t\geq 0\}$, 
$N^{\lambda}=\{N_t^{\lambda}\colon t\geq 0\}$
and 
$N^{\mu}=\{N_t^{\mu}\colon t\geq 0\}$
two Poisson processes with intensities $\lambda$ and $\mu$,
respectively, 
$Y^{\alpha}=\{Y_i^{\alpha}:i=1,2,\dots\}$ 
and 
$Y^{\beta}=\{Y_i^{\beta}:i=1,2,\dots\}$
two sequences of independent exponentially distributed random variables
with parameters $\alpha$ and $\beta$, respectively. 
Moreover, $W,N^{\lambda},N^{\mu},Y^{\alpha}$ and $Y^{\beta}$
are assumed to be independent.
The process $X=\{X_t\}$ is now defined via
\begin{equation}
\label{e:process}
X_t=at+b W_t+\sum_{i=1}^{N_t^{\lambda}}Y_i^{\alpha}-\sum_{i=1}^{N_t^{\mu}}Y_i^{\beta},
\end{equation}
where $a$ and $b\geq 0$ are real parameters. Clearly, $X$ is a L\'evy process
and its L\'evy-Khintchine representation is given by
\begin{equation}
\label{e:lk}
\E\left(\exp (z\, X_t)\right)=\exp (t\,\psi(z))
\end{equation}
with 
\begin{equation}
\label{e:psi}
\psi(z)=az+\frac 12b^2z^2+\lambda\,\frac{z}{\alpha-z}-\mu\,\frac{z}{\beta+z}.
\end{equation}
It is enough for our purposes to take hereby $z$ real, and then the
representation in (\ref{e:lk}) holds for $z\in(-\beta,\alpha)$.

Next we compute the Green kernel of $X$ when 
all the parameters in \eqref{e:psi} are positive.
It is easily seen that for $r\geq 0$
the equation $\psi(z)=r$ has exactly four solutions $\rho_k,\,
k=1,2,3,4.$ These satisfy
\begin{equation}
\label{roots} 
\rho_1<-\beta<\rho_2\leq 0<\rho_3<\alpha<\rho_4
\end{equation}
and
\begin{equation}
\label{roots2} 
\psi'(\rho_1)<0,\ \ \psi'(\rho_2)<0,\ \ \psi'(\rho_3)>0,\ \ \psi'(\rho_4)>0.
\end{equation}
Notice that $\rho_2=0$ if and only if $r=0,$ in which case it is
assumed $\psi'(0)<0$ implying 
$$
\lim_{t\to\infty}X_t=-\infty\quad {\rm a.s.}
$$
Using the general
definition of the resolvent, see  (\ref{resolvent}), 
we have for $z\in(\rho_2,\rho_3)$ 
\begin{eqnarray}
\label{e:expansion}
&&\hskip-1cm
\nonumber
\int_{-\infty}^\infty {\rm e}^{z\,x}\,G_r(0,dx)
=
\int_0^\infty dt\,{\rm e}^{-rt}\,\E\left(\exp(z\, X_t)\right)
=\frac
1{r-\psi(z)}\\
&&\hskip2.3cm
\nonumber
=
\frac{\psi'(\rho_1)^{-1}}{\rho_1-z}+
\frac{\psi'(\rho_2)^{-1}}{\rho_2-z}+
\frac{\psi'(\rho_3)^{-1}}{\rho_3-z}+
\frac{\psi'(\rho_4)^{-1}}{\rho_4-z}.
\end{eqnarray}
Consequently, inverting the right hand side yields 
\begin{equation}
\label{e12}
G_r(0,dx)=
\begin{cases}
-{\psi'(\rho_1)^{-1}}{\rm e}^{\,-\rho_1x}dx-
{\psi'(\rho_2)^{-1}}{\rm e}^{\,-\rho_2x}dx,& x< 0,\\ 
{\psi'(\rho_3)^{-1}}{\rm e}^{\,-\rho_3x}dx+
{\psi'(\rho_4)^{-1}}{\rm e}^{\,-\rho_4x}dx,& x> 0.\\ 
\end{cases}
\end{equation}
and, hence, the resolvent is absolutely continuous with respect
to Lebesgue measure. With slight abuse of notation, we let $G_r(0,x)$ denote also 
the Green kernel, i.e., the density of the resolvent $G_r$ with
respect to the Lebesgue measure. From the spatial homogeniety of $X$
it follows that $G_r(x,0)=G_r(0,-x).$

The absolute continuity of the resolvent can
alternatively  be verified by checking that the condition (ii) in Theorem II.5.16 in
Bertoin \cite{bertoin96} 
holds . We recall also the general result (see \cite{bertoin96} p. 25) 
which says that the absolute continuity of the resolvent is equivalent with the
property that $x\mapsto G_rf(x)$ is continuous for all essentially bounded
measurable functions $f.$ 

As we have noticed, the process
$\widehat{X}=\{-X_t\}$ may be viewed as a dual process 
associated with $X.$ Let $\widehat G_r$ denote the resolvent of
$\widehat{X}$. Then the duality relationship
$$
\int dx\, f(x)\, G_r g(x)=\int dx\, \widehat G_r f(x)\,g(x)
$$  
holds the duality measure being the Lebesgue measure. 
The Green kernel  of the dual process  
is given by
$$
\widehat G_r(x,y)=G_r(y,x).
$$
Notice that the value of 
$x\mapsto G_r(0,x)$ at 0 is chosen so that the
resulting function is lower 
semi-continuous (since the Green kernel when considered as a function
of the second argument should be excessive for the dual process).

To conclude the above discussion, we have verified Hypothesis (B) in
\cite{kunitawatanabe65}, that is, ${\rm (h_1)}, {\rm (h_2)}$ and ${\rm
  (h_3)}$ in Section 2 are fullfilled. Consequently, also  ${\rm (p_1)}, {\rm (p_2)}$ 
and ${\rm  (p_3)}$ in Section 2 are valid and the Riesz decomposition 
(\ref{riesz}) holds. Moreover, it can be proved, e.g. 
using the Martin boundary theory, as presented in \cite{kunitawatanabe65}, that the harmonic
function $h_r$  appearing in (\ref{riesz}) is of the form      
$$
h_r(x)=c_1\,{\rm e}^{\,\rho_2x}+ c_2\,{\rm e}^{\,\rho_3x},  
$$
where  $c_1$ and $c_2$ are non-negative constants.

It is interesting to note that 
when multiplying both sides of
\eqref{e:expansion}
by $z$ and letting $z\to\infty$
we obtain, in case $b>0$ (cf. (\ref{e:psi})), 
\begin{equation*}
\frac 1{\psi'(\rho_1)}+ 
\frac 1{\psi'(\rho_2)}+
\frac 1{\psi'(\rho_3)}+
\frac 1{\psi'(\rho_4)}=0,
\end{equation*}
which implies that the Green kernel
is continuous at $x=0$.
But, when $b=0$, the Green kernel may be 
discontinuous. This happens, for instance, when 
$X$ is a compound Poisson process with negative drift and
exponentially distributed positive jumps.
More precisely, taking $b=\mu=0$, and $a<0$ in \eqref{e:process}  
the characteristic exponent reduces to
$$
\psi(z)=az+\lambda\,\frac{z}{\alpha-z}.
$$
Now there are only two roots 
$\rho_1$ and $\rho_2$
and these satisfy
$$
\rho_1\leq 0< \rho_2.
$$
Consequently, 
\begin{equation}
\label{e:psi2}
\frac
1{r-\psi(z)}
=
\frac{\psi'(\rho_1)^{-1}}{\rho_1-z}+
\frac{\psi'(\rho_2)^{-1}}{\rho_2-z},
\end{equation}
and we have the Green kernel
\begin{equation}
\label{e121}
G_r(0,x)=
\begin{cases}
-{\psi'(\rho_1)^{-1}}\,{\rm e}^{\,-\rho_1x},& x< 0,\\ 
{\psi'(\rho_2)^{-1}}{\rm e}^{\,-\rho_2x},& x\geq 0.\\ 
\end{cases}
\end{equation}
From \eqref{e:psi2} it is seen that 
$$
\psi'(\rho_1)^{-1}+\psi'(\rho_2)^{-1}=\frac 1a,
$$ 
and, hence, $x\mapsto G_r(0,x)$ is discontinuous at 0 (but lower
semi-continuous since $a<0$ implies $-\psi'(\rho_1)<\psi'(\rho_2)$).

\subsection{Optimal stopping of processes with two sided exponential Green kernel}
\label{sec42}

We consider here a subclass of processes introduced in Section
\ref{sec41} the aim being to apply results in Theorem \ref{prop1}
and \ref{prop2}. Indeed, let $X$ be a L\'evy process having a
Green kernel with the following simple exponential structure:  
\begin{equation}
\label{g1}
G_r(x):=G_r(0,x)=
\begin{cases}
-A_1\,{\rm e}^{-\rho_1\, x},&\text{ $x<0$},\\\
A_2\, {\rm e}^{-\rho_2\, x},&\text{ $x\geq0$},\\
\end{cases}
\end{equation}
where $\rho_{1,2}$ are the roots of
the equation $\psi(z)=r$ such that $\rho_1\leq 0<\rho_2$ and 
$A_{1,2}=1/\psi'(\rho_{1,2}).$ 
In the case $r=0$ it is assumed that the process drifts to $-\infty$
and, hence, we have $\rho_1=0.$

The Green kernel of form (\ref{g1}) appears in two basic cases which, 
using the notation in (\ref{e:lk}) and (\ref{e:psi}), are:
\begin{description}
\item{$\bullet$}\hskip5mm Wiener process with drift, i.e., $b>0,$ $\lambda=\mu=0,$
 
\item{$\bullet$}\hskip5mm compound Poisson process with negative drift and positive
exponential jumps, i.e., $a<0,$ $b=0,$ $\lambda>0$ and $\mu=0.$ 
\end{description}
The point we want to make here is that our approach to optimal
stopping treats these processes similarly. Recall that in the case of
Wiener process usually smooth pasting is valid when moving from the
continuation region to the stopping region but in the compound
Poisson case there is ``only'' continuous pasting. In other words, our
approach does not use smooth pasting as a tool,
but this property can, of course, be checked (when valid) 
from the calculated explicit form
of the value function.

\begin{proposition}
\label{prop51}
For a given $x^*\in\R$ let $\sigma$ be a measure on $[x^*,+\infty)$
  with a continuously differentiable density $\sigma'$ on $(x^*,+\infty).$ Then 
the function  
\begin{equation*}
V(x):=\int_{x^*}^{\infty}G_r(y-x)\sigma(dy)
\end{equation*}
is two times continuously differentiable on
$D:=\{x\in\R\colon x\neq x^*\}$ 
and satisfies on $D$ the ordinary differential equation (ODE) 
\begin{eqnarray}
\label{sigma1}
&&\hskip-1cm
V''(x)-(\rho_2+\rho_1)V'(x)+\rho_1\rho_2 V(x)
\\
&&\hskip2cm
\nonumber
=-(A_2+A_1)\sigma''(x)+(\rho_2 A_1+\rho_1A_2)\sigma'(x),
\end{eqnarray}
where $\sigma''(x)=\sigma'(x)=0$ for $x\in(-\infty,x^*).$
\end{proposition}
\begin{proof}
From the definition of $V$, taking into account the
form of the Green kernel, we have for $x>x^*$
\begin{equation*}
V(x)=-A_1 {\rm e}^{\rho_1 x}\int_{x^*}^{x}{\rm e}^{-\rho_1 y}\sigma(dy)
+A_2 {\rm e}^{\rho_2 x}\int_{x}^{\infty}{\rm e}^{-\rho_2 y}\sigma(dy).
\end{equation*}
The right hand side of this equation can be differentiated twice
proving that $V''$ exist in $D,$ and the claimed ODE is obtained after
some straightforward manipulations.
\end{proof}

\begin{corollary}
\label{cor51}
Let $X$ be a Wiener process with drift. Then the ODE in (\ref{sigma1})
takes the form 
\begin{equation}
\label{sigma3}
\frac{b^2}2\,V''(x)+a\,V'(x)-r\, V(x)=-(a^2+2b^2r)\,\sigma'(x).
\end{equation}
\end{corollary}
\begin{proof}
The quantities needed to derive (\ref{sigma3}) from   (\ref{sigma1}) are
  $$
\rho_1=-\frac 1{b^2}\left(\sqrt{a^2+2b^2 r}+a\right),\quad
\rho_2=\frac 1{b^2}\left(\sqrt{a^2+2b^2 r}-a\right)
$$
and
$$
A_1=-\sqrt{a^2+2b^2 r},\quad A_2=\sqrt{a^2+2b^2 r}.
$$
In particular, 
notice that  $A_2+A_1=0$ which reflects the fact that the Green kernel is continuous.
\end{proof}
In Novikov and Shiryayev \cite{novikovshiryaev04}
the optimal stopping problem for a
general random walk with reward function $\max\{0, x^n\}$,
$n=1,2,\dots$, is considered, and the solution is characterized via the
Appell polynomials associated with the distribution of the
maximum of the process. 
In the next example we present explicit reults for a more
general reward function, that is,  $\max\{0, x^\gamma\}$, $\gamma\geq 1,$ 
but for a more particular L\'evy process studied in the subsection. 

\begin{example}
\label{ex2} {\rm
Let $X$ denote a compound Poisson process with negative drift and positive
exponential jumps, i.e., take $a<0,$ $b=0$ and $\mu=0$ in
(\ref{e:process}). 
For simplicity, we consider optimal stopping problem without
discounting: 
$$
\sup_{\tau\in{\cal M}}\E_x\left(g(X_{\tau})\right),
$$
where $g(x):=\max\{0, x^\gamma\}$ with $\gamma\geq 1.$ For $r=0$ the Green kernel of
$X$ is
\begin{equation}
\label{e13}
G(x,0):=G_0(x,0)=
\begin{cases}
A_2\,{\rm e}^{\rho\,x}\,dx,& x\leq 0,\\ 
-A_1,& x> 0,\\ 
\end{cases}
\end{equation}
where
\begin{equation}
\label{e131}
 \rho:=\rho_2=\al+\frac\lambda a>0
\end{equation}
and
$$
A_1=\frac\al{\la+a\al}<0,\qquad A_2=\frac\lambda{a(
  \lambda+a\al)}>0.
$$
Notice that $\rho>0$ means that a.s. $\lim_{t\to\infty}X_t=-\infty.$ 

Our aim is to find a measure $\sigma$ and a number $x^*$ such that the function $V$
defined via 
\begin{equation}
\label{v}
V(x)=\int_{[x^*,+\infty)}G(x,y)\sigma(dy)  
\end{equation}
has properties (a), (b), (c) and (d) given in Theorem \ref{prop1}.

To begin with,  
consider equation (\ref{sigma1}) for $x>x^*$ and
$V(x)=x^\gamma,$  that is,
$$
-\sigma''(x)+\alpha\sigma'(x)=a\gamma(\gamma-1)x^{\gamma-2}-\left(a\alpha+\lambda\right)
\gamma x^{\gamma-1}.
$$     
Assuming $\lim_{x\to+\infty} {\rm e}^{-\alpha x}\,\sigma'(x)=0$ we
obtain the solution 
$$
\sigma'(x)=-a\gamma\,x^{\gamma-1}-\lambda {\rm e}^{\alpha
  x}\int_x^\infty{\rm e}^{-\alpha y} \gamma\,y^{\gamma-1}\, dy.
$$
If $\gamma =1$ then $\sigma'(x)=-a-(\lambda/\alpha)>0.$ For $\gamma>1$
it is easily seen that $\sigma'(0)<0$ and $\sigma'(x)\to+\infty$ as $x\to\infty.$

The claim is that the equation $\sigma'(x)=0,$ that is 
\begin{equation}
\label{xo}
x^{\gamma-1}=\frac\lambda{(-a)}\,
 {\rm e}^{\alpha x}
\int_x^\infty {\rm e}^{-\alpha z}\,z^{\gamma-1}\, dz,
\end{equation}
has a unique solution for $x>0,$ which we denote by $x^*_{\gamma-1}.$  
Equation (\ref{xo}) is equivalent to
\begin{equation}
\label{eq:xg}
F(x,\gamma-1)=1,
\end{equation}
if we define
\begin{equation}\label{eq:ef}
F(x;{u}):=\frac\lambda{(-a)}\,\int_0^\infty{\rm e}^{-\alpha y} \,\left(1+\frac yx\right)^{u}\,dy.
\end{equation}
We revise some properties of the function just introduced.
\begin{lemma}
The function $F(x,{u})$ in \eqref{eq:ef} defines an implicit function $\varphi\colon[1,\infty)\to\R$ 
such that $F(\varphi(u),{u})=1$ for each $u\geq 1$.
Furthermore, the function $\varphi$ is increasing, and
satisfies the inequality
\begin{equation}\label{eq:inequalities}
\varphi(1)<\varphi(u)<\frac u\rho.
\end{equation}
\end{lemma}
\begin{proof}
It is not difficult to verify that, 
for fixed ${u}>0$,
the function $F$ is decreasing in $x$, and that
$$
\lim_{x\to 0+}F(x,u)=\infty,
\qquad
\lim_{x\to\infty}F(x,u)=\frac\lambda{(-a)\alpha}<1.
$$
This means that for any $u\geq 1$ the equation $F(x,{u})=1$
has a unique solution $x:=\varphi(u)$. 
Furthermore, it is also clear that,
for fixed $x>0$, the function $F(x,u)$ is increasing in ${u}$.
This means, that $\varphi$
is increasing, as
$$
\frac{\partial\varphi}{\partial u}=
-\frac{\partial F}{\partial x}
\left(\frac{\partial F}{\partial u}\right)^{-1}>0.
$$
Finally, multiplying  the inequality 
$$
\left(1+u\frac yx\right)\leq
\left(1+\frac yx\right)^{u}\leq 
{\rm e}^{uy/x}
$$
by ${\rm e}^{-\alpha y}$ and integrating we obtain that 
$$
F_1(x,u):=
\frac\lambda{(-a)\alpha}
\left(1+\frac u{\alpha x}\right)
<
F(x,u)
<
\frac\lambda{(-a)}\frac 1{\alpha-u/x}=:F_2(x,u).
$$
and \eqref{eq:inequalities}
follows as the bounds are the respective roots of
the equations $F_1(x,u)=1$, $F_2(x,u)=1$,
and, in particular, the root $x_1$ of the first equation is
$$
x_1=\frac{\gamma\lambda}{(-a)\alpha\rho},
$$
and $x^*_1=-\lambda/a\alpha\rho$. This last value
can be computed from the equation $F(x,1)=1$, 
and was found in  \cite{mordecki99}.
This concludes the proof of the Lemma.
\end{proof}

Observe now that for $x>x^*_{\gamma-1}=:x^o_\gamma$
the function $\sigma'$ induces a positive Radon
measure on $(x^o_\gamma,\infty)$. 
However, since, for any constant $c,$ the
function $x^\gamma+c$ induces the same measure as the function
$x^\gamma$ it remains to find, 
the support of $\sigma$ of the form $(x^*,\infty)$ such that for all
$x>x^*$ 
$$
x^\gamma =\int_{[x^*,+\infty)}G(x,y)\,\sigma'(y)\,dy.  
$$
Therefore, consider
\begin{eqnarray*}
&&
\hskip-.8cm
\int_x^\infty G(x,y)\sigma(dy)
=A_2\,  {\rm e}^{\rho x}    
\int_x^\infty {\rm e}^{-\rho y}\,\sigma'(y)\,dy\\
&&
\hskip1.5cm
=
\gamma\,A_2\,  {\rm e}^{\rho x}    
\int_x^\infty {\rm e}^{-\rho y}\,\left(
-a\,y^{\gamma-1}-\lambda {\rm e}^{\alpha
  y}\int_y^\infty{\rm e}^{-\alpha z}\,z^{\gamma-1}dz\right)\, dy.
\end{eqnarray*}
Applying Fubini's theorem for the latter term yields
\begin{eqnarray*}
&&
\hskip-.8cm
\int_x^\infty dy\,{\rm e}^{(\alpha-\rho) y}\,
\int_y^\infty dz\, {\rm e}^{-\alpha z}\,z^{\gamma-1}
\\
&&\hskip2cm
=
\frac 1{\alpha-\rho}\left(
\int_x^\infty {\rm e}^{-\rho z}\,z^{\gamma-1}\, dz
-
 {\rm e}^{(\alpha-\rho)x}
\int_x^\infty  {\rm e}^{-\alpha z}\,z^{\gamma-1}\, dz
\right).
\end{eqnarray*}
Observing that $-a=\lambda/(\alpha-\rho)$ we have 
\begin{eqnarray*}
&&
\hskip-1cm
\int_x^\infty G(x,y)\sigma(dy)
=\frac {\alpha-\rho}{\rho}\,
 {\rm e}^{\alpha x}
\int_x^\infty {\rm e}^{-\alpha z}\,\gamma\,z^{\gamma-1}\, dz
\end{eqnarray*}
Consequently, after an integration by parts, the equation
$$
x^\gamma=\frac {\alpha-\rho}{\rho}\,
 {\rm e}^{\alpha x}
\int_x^\infty {\rm e}^{-\alpha z}\,\gamma\,z^{\gamma-1}\, dz
$$
is seen to be equivalent with 
\begin{equation}
\label{x*}
x^\gamma={(\alpha-\rho)}\,
 {\rm e}^{\alpha x}
\int_x^\infty {\rm e}^{-\alpha z}\,z^{\gamma}\, dz
\end{equation}
which coincides with equation (\ref{xo}) if therein $\gamma-1$ is
changed to $\gamma.$ Hence, equation (\ref{x*}) has a unique solution
which is, using the notation introduced above, $x^*_\gamma.$ As the function $x=\varphi({u})$ 
is increasing, we know that $x^o_\gamma=\varphi(\gamma-1)<\varphi(\gamma)=x^*_\gamma$.

Next step is to verify that the value function obtained from
(\ref{v}) is continuous and satisfies $V(x)>x^\gamma$ for
$x<x^*_\gamma.$ Therefore consider for $x<x^*_\gamma$ 
\begin{equation}
\label{v1}
V(x)=\int_{[x^*_\gamma,+\infty)}G(x,y)\sigma(dy)  
={\rm e}^{\rho(x-x^*_\gamma)}(x^*_\gamma )^\gamma.
\end{equation}
Consequently, $V$ is continuous and 
\begin{equation}
\label{majo}
V(x)>x^\gamma\quad \Leftrightarrow \quad {\rm e}^{-\rho x^*_\gamma}(x^*_\gamma
)^\gamma> {\rm e}^{-\rho x}x^\gamma.
\end{equation}
The right hand side of (\ref{majo}) holds if $x\mapsto G(x):= {\rm e}^{-\rho x}x^\gamma$ is
increasing for $x<x^*_\gamma.$ Clearly, $G'$ is positive if 
$$
-\rho x+\gamma >0,
$$
and this holds since $x^*_\gamma<\gamma/\rho$ by the second inequality in \eqref{eq:inequalities}.

To conclude, the optimal stopping time
$\tau^*$ is given by
$$
\tau^*:=\inf\{t : X_t\geq x^*_\gamma\},
$$
and the value function $V$ is 
for $x<x^*_\gamma$ as in (\ref{v1}). Since $V'(x^*_\gamma-)=\rho (x^*_\gamma)^\gamma$
and $g'(x^*_\gamma)=\gamma (x^*_\gamma)^{\gamma-1}$ there is no smooth fit at $x^*_\gamma$.

We conclude by presenting the following table with some numerical
results. The computations are done with Mathematica-package where one
can find a subroutine for incomplete gamma-function and programs for
numerical solutions of equations based on standard Newton-Raphson's method and  
the secant method. A good starting value for Newton-Raphson's method
seems to be $\gamma/\rho.$ It is interesting to notice from the table that if $\rho<<\alpha$
then $x^*_\gamma\simeq \gamma/\rho.$

$$
\vbox{\offinterlineskip\halign{&\vrule#&\strut\quad\hfil
$#$\quad\hfil\cr
\noalign{\hrule}
height10pt&\omit&&\omit&&\omit&&\omit&&\omit&&\omit&&\omit&\cr
&{\rm \alpha}&&{\rm \rho}&&{-\lambda/a}&&{\gamma}&&{\gamma/\rho}&&{x^*_\gamma}&&{x^o_\gamma}  &\cr
height10pt&\omit&&\omit&&\omit&&\omit&&\omit&&\omit&&\omit&\cr
\noalign{\hrule height0.1pt}
height10pt&\omit&&\omit&&\omit&&\omit&&\omit&&\omit&&\omit&\cr
&
10 
&&
1
&&
9
&&
20
&&
20
&& 
19.8896
&&
18.8896
&
\cr
height10pt&\omit&&\omit&&\omit&&\omit&&\omit&&\omit&&\omit&\cr
&
10 
&&
1
&&
9
&&
10
&&
10
&& 
9.8902
&&
8.8904
&
\cr
height10pt&\omit&&\omit&&\omit&&\omit&&\omit&&\omit&&\omit&\cr
&
10 
&&
1
&&
9
&&
5
&&
5
&& 
4.8915
&&
3.8921
&
\cr
height10pt&\omit&&\omit&&\omit&&\omit&&\omit&&\omit&&\omit&\cr
&
10 
&&
1
&&
9
&&
2.5
&&
2.5
&& 
2.3939
&&
1.3968
&
\cr
height10pt&\omit&&\omit&&\omit&&\omit&&\omit&&\omit&&\omit&\cr
&
10 
&&
1
&&
9
&&
1
&&
1
&& 
.9
&&
-
&
\cr
height10pt&\omit&&\omit&&\omit&&\omit&&\omit&&\omit&&\omit&\cr
\noalign{\hrule height0.1pt}
height10pt&\omit&&\omit&&\omit&&\omit&&\omit&&\omit&&\omit&\cr
&
10 
&&
9
&&
1
&&
20
&&
2.2222
&& 
1.7613
&&
1.6579
&
\cr
height10pt&\omit&&\omit&&\omit&&\omit&&\omit&&\omit&&\omit&\cr
&
10 
&&
9
&&
1
&&
10
&&
1.1111
&& 
.7511
&&
.6547
&
\cr
height10pt&\omit&&\omit&&\omit&&\omit&&\omit&&\omit&&\omit&\cr
&
10 
&&
9
&&
1
&&
5
&&
.5555
&& 
.2881
&&
.2045
&
\cr
height10pt&\omit&&\omit&&\omit&&\omit&&\omit&&\omit&&\omit&\cr
&
10 
&&
9
&&
1
&&
2.5
&&
.2789
&& 
.0917
&&
.0319
&
\cr
height10pt&\omit&&\omit&&\omit&&\omit&&\omit&&\omit&&\omit&\cr
&
10 
&&
9
&&
1
&&
1
&&
.1111
&& 
.0111
&&
-
&
\cr
height10pt&\omit&&\omit&&\omit&&\omit&&\omit&&\omit&&\omit&\cr
\noalign{\hrule height0.1pt}
height10pt&\omit&&\omit&&\omit&&\omit&&\omit&&\omit&&\omit&\cr
&
1 
&&
.5
&&
.5
&&
20
&&
40
&& 
38.1592
&&
36.166
&
\cr
height10pt&\omit&&\omit&&\omit&&\omit&&\omit&&\omit&&\omit&\cr
&
1 
&&
.5
&&
.5
&&
10
&&
20
&& 
18.2726
&&
16.2942
&
\cr
height10pt&\omit&&\omit&&\omit&&\omit&&\omit&&\omit&&\omit&\cr
&
1 
&&
.5
&&
.5
&&
5
&&
10
&& 
8.4369
&&
6.5011
&
\cr
height10pt&\omit&&\omit&&\omit&&\omit&&\omit&&\omit&&\omit&\cr
&
1 
&&
.5
&&
.5
&&
2.5
&&
5
&& 
3.6529
&&
1.8398
&
\cr
height10pt&\omit&&\omit&&\omit&&\omit&&\omit&&\omit&&\omit&\cr
&
1 
&&
.5
&&
.5
&&
1
&&
2
&& 
1
&&
-
&
\cr
height10pt&\omit&&\omit&&\omit&&\omit&&\omit&&\omit&&\omit&\cr
\noalign{\hrule}}}
$$
}
\end{example}

\subsection*{Acknowledgements}
The second author thanks \AA bo Akademi (\AA bo, Finland),
and the National Visitors Program (Finland) for the support and hospitality,
that made possible the initiation of this work.
\bibliographystyle{plain}
\bibliography{pape1}

\begin{thebibliography}{10}

\bibitem{bertoin96}
J.~Bertoin.
\newblock {\em L\'evy Processes}.
\newblock Cambridge University Press, Cambridge, 1996.

\bibitem{blumenthalgetoor68}
R.M. Blumenthal and R.K. Getoor.
\newblock {\em Markov {P}rocesses and {P}otential {T}heory}.
\newblock Academic Press, New York, London, 1968.

\bibitem{boyarchenkolevendorskij02}
S.~I. Boyarchenko and S.~Z. Levendorskij.
\newblock {\em Non-Gaussian Merton-Black-Scholes theory}.
\newblock World Scientific, Singapore, 2002.

\bibitem{chan00}
T.~Chan.
\newblock Pricing perpetual {A}merican options driven by spectrally one-sided
  {L}\'evy processes.
\newblock {\em Preprint}, 2000.

\bibitem{darlingliggettaylor72}
D.A. Darling, T.~Liggett, and H.M. Taylor.
\newblock {Optimal stopping for partial sums.}
\newblock {\em Ann. Math. Stat.}, 43:1363--1368, 1972.

\bibitem{dynkin63}
E.B. Dynkin.
\newblock The optimum choice of the instant for stopping a markov process.
\newblock {\em Doklady {A}kademii {N}auk {SSSR}}, 150(2):238--240, 1963.

\bibitem{dynkin69}
E.B. Dynkin.
\newblock Prostranstvo vyhodov markovskogo processa.
\newblock {\em Uspehi Mat. Nauk (English trans. in Russ. Math. Surv.)}, XXIV 4
  (148):89--152, 1969.

\bibitem{folland99}
G.B. Folland.
\newblock {\em Real Analysis}.
\newblock Wiley, New York, 2nd. edition, 1999.

\bibitem{gerbershiu98}
H.U. Gerber and E.S.W. Shiu.
\newblock Pricing perpetual options for jump processes.
\newblock {\em North Amer. Act. J.}, 2(3):101--112, 1998.

\bibitem{kouwang04}
S.G. Kou and H.T. Wang.
\newblock Option pricing under a double exponential jump diffusion model.
\newblock {\em Manag. Sci.}, 50(9):1178--1192, 2004.

\bibitem{kunitawatanabe65}
H.~Kunita and T.~Watanabe.
\newblock Markov processes and {M}artin boundaries,~{I}.
\newblock {\em Illinois J. Math.}, 9(3):485--526, 1965.

\bibitem{kyprianousurya05}
A.~E. Kyprianou and B.~A. Surya.
\newblock On the {N}ovikov-{S}hiryaev optimal stopping problems in continous
  time.
\newblock {\em Electronic Communications in Probability}, 10:146--154, 2005.

\bibitem{mckean67}
H.P. McKean.
\newblock A free boundary problem for the heat equation arising from a problem
  in mathematical economics.
\newblock {\em Indust. Management Rev}, 6:32--39, 1965.

\bibitem{mordecki99}
E.~Mordecki.
\newblock Optimal stopping for a compound poisson process with exponential
  jumps.
\newblock {\em Publicaciones Matem\'aticas del Uruguay}, 7:55--66, 1997.

\bibitem{mordecki02}
E.~Mordecki.
\newblock {Optimal stopping and perpetual options for {L}\'evy processes.}
\newblock {\em Finance Stoch.}, 6(4):473--493, 2002.

\bibitem{mordecki02a}
E.~Mordecki.
\newblock Perpetual options for {L}\'evy processes in the {B}achelier model.
\newblock In A.N. Shiryaev, editor, {\em Proc. Steklov Inst. Math. Stochastic
  Financial Mathematics}, volume 237, pages 247--255. Nauka, Moscow, 2002.

\bibitem{novikov06}
A.~Novikov and Shiryaev A.
\newblock On the solution of the optimal stopping problem for random walks and
  {L}\'evy process.
\newblock Conference at the Optimal Stopping with Applications Symposium,
  Manchester, UK, January 2006.

\bibitem{novikovshiryaev04}
A.~Novikov and A.N. Shiryaev.
\newblock On an effective solution of the optimal stopping problem for random
  walks.
\newblock {\em Th. Probab. Appl.}, 49:373--382, 2004.

\bibitem{rogozin66}
B.A. Rogozin.
\newblock {On distributions of functionals related to boundary problems for
  processes with independent increments}.
\newblock {\em Th. Probab. Appl.}, 11:580--591, 1966.

\bibitem{salminen85}
P.~Salminen.
\newblock Optimal stopping of one-dimensional diffusions.
\newblock {\em Math.\ Nachr.}, 124:85--101, 1985.

\bibitem{sato99}
K.~Sato.
\newblock {\em L\'evy processes and infinitely divisible distributions}.
\newblock Cambridge Press, Cambridge, 1999.

\bibitem{shiryayev78}
A.~N. Shiryayev.
\newblock {\em Optimal {S}topping {R}ules}.
\newblock Springer {V}erlag, New York, Heidelberg, Berlin, 1978.

\bibitem{snell53}
L.~Snell.
\newblock Applications of martingale system theorems.
\newblock {\em Trans. Amer. Math. Soc.}, 73:293--312, 1953.

\end{thebibliography}

\end{document}